\documentclass[a4paper,reqno,12pt]{amsart}
\usepackage{amsmath,amsrefs}

\numberwithin{equation}{section}

\DeclareMathOperator{\divergence}{div}
\DeclareMathOperator{\dist}{dist}
\DeclareMathOperator{\meas}{meas}
\DeclareMathOperator{\loc}{loc}
\DeclareMathOperator{\supp}{supp}
\DeclareMathOperator{\sgn}{sgn}
\newcommand{\OO}{\operatorname{O}}
\newcommand{\oo}{\operatorname{o}}
\newcommand{\R}{\mathbb{R}}
\newcommand{\N}{\mathbb{N}}

\renewcommand{\[}{\left[}
\renewcommand{\]}{\right]}
\renewcommand{\(}{\left(}
\renewcommand{\)}{\right)}

\newtheorem{theorem}{Theorem}[section]
\newtheorem{corollary}[theorem]{Corollary}
\newtheorem{lemma}[theorem]{Lemma}
\newtheorem{remark}[theorem]{Remark}

\begin{document}

\date{July 23, 2014}

\title[A priori estimates for critical $p$--Laplace equations]{A priori estimates and application to the symmetry of solutions for critical $p$--Laplace equations}

\author{J\'er\^ome V\'etois}
\address{J\'er\^ome V\'etois, Univ. Nice Sophia Antipolis, CNRS,  LJAD, UMR 7351, 06108 Nice, France.}
\email{vetois@unice.fr}

\begin{abstract}
We establish pointwise a priori estimates for solutions in $D^{1,p}\(\R^n\)$ of equations of type $-\Delta_pu=f\(x,u\)$, where $p\in\(1,n\)$, $\Delta_p:=\divergence\big(\left|\nabla u\right|^{p-2}\nabla u\big)$ is the $p$--Laplace operator, and $f$ is a Caratheodory function with critical Sobolev growth. In the case of positive solutions, our estimates allow us to extend previous radial symmetry results. In particular, by combining our results and a result of  Damascelli--Ramaswamy~\cite{DamRam}, we are able to extend a recent result of Damascelli--Merch\'an--Montoro--Sciunzi~\cite{DamMerMonSci} on the symmetry of positive solutions in $D^{1,p}\(\R^n\)$ of the equation $-\Delta_pu=u^{p^*-1}$, where $p^*:=np/\(n-p\)$.
\end{abstract}

\maketitle

\section{Introduction and main results}\label{Sec1}

In this paper, we are interested in problems of the type
\begin{equation}\label{Eq1}
\left\{\begin{aligned}&-\Delta_pu=f\(x,u\)\quad\text{in }\R^n,\\
&u\in D^{1,p}\(\R^n\),
\end{aligned}\right.
\end{equation}
where $p\in\(1,n\)$, $\Delta_pu:=\divergence\(\left|\nabla u\right|^{p-2}\nabla u\)$, $D^{1,p}\(\R^n\)$ is the completion of $C^\infty_c\(\R^n\)$ with respect to the norm $\left\|u\right\|_{D^{1,p}\(\R^n\)}:=\(\int_{R^n}\left|\nabla u\right|^pdx\)^{1/p}$, and $f:\R^n\times\R\to\R$ is a Caratheodory function such that
\begin{equation}\label{Eq2}
\left|f\(x,s\)\right|\le\Lambda\left|s\right|^{p^*-1}\quad\text{for all }s\in\R\text{ and a.e. }x\in\R^n, 
\end{equation}
for some real number $\Lambda>0$, with $p^*:=np/\(n-p\)$.

\smallskip
Our main result is as follows.

\begin{theorem}\label{Th}
Let $p\in\(1,n\)$, $f:\R^n\times\R\to\R$ be a Caratheodory function such that \eqref{Eq2} holds true and $u$ be a solution of \eqref{Eq1}. Then there exists a constant $C_0=C_0\(n,p,\Lambda,u\)$ such that 
\begin{equation}\label{ThEq1}
\left|u\(x\)\right|\le C_0\big(1+\left|x\right|^{\frac{n-p}{p-1}}\big)^{-1}\quad\text{and}\quad\left|\nabla u\(x\)\right|\le C_0\big(1+\left|x\right|^{\frac{n-1}{p-1}}\big)^{-1}
\end{equation} 
for all $x\in\R^n$. If moreover $u\ge0$ in $\R^n$ and $\int_{\R^n}f\(x,u\)dx>0$, then we have
\begin{equation}\label{ThEq2}
u\(x\)\ge C_1\big(1+\left|x\right|^{\frac{n-p}{p-1}}\big)^{-1}
\end{equation}
for all $x\in\R^n$, for some constant $C_1=C_1\(n,p,\lambda,\Lambda,u\)>0$, where $\lambda$ is a real number such that $0<\lambda<\int_{\R^n}f\(x,u\)dx$.
\end{theorem}

The dependence on $u$ of the constants $C_0$ and $C_1$ will be made more precise in Remarks~\ref{Rem1} and~\ref{Rem2}.

\smallskip
In the case of the Laplace operator ($p=2$), the upper bound estimates \eqref{ThEq1} have been established by Jannelli--Solimini~\cite{JanSol} for nonlinearities of the form $f\(x,u\)=\sum_{i=1}^Na_i\(x\)\left|u\right|^{q_i^*-2}u$, where $q_i^*:=2^*\(1-1/q_i\)$, $q_i\in\(n/2,\infty\]$, $\left|a_i\(x\)\right|=\OO\big(\left|x\right|^{-n/q_i}\big)$ for large $\left|x\right|$, and $a_i$ belongs to the Marcinkiewicz space $M^{q_i}\(\R^n\)$ for all $i=1,\dotsc,N$. The case of unbounded domains $\Omega\ne\R^n$ is also treated in~\cite{JanSol}. 

\smallskip
Since the pioneer work of Gidas--Ni--Nirenberg~\cite{GidNiNir} and later extensions by Li~\cite{Li} in case $p=2$ and Damascelli--Ramaswamy~\cite{DamRam} in case $1<p<2$, decay estimates are known to be useful to derive radial symmetry results for $C^1$--solutions of problems of the type 
\begin{equation}\label{Eq3}
\left\{\begin{aligned}&-\Delta_pu=f\(u\),\quad u>0\quad\text{in }\R^n,\\
&u\(x\)\longrightarrow0\quad\text{as }\left|x\right|\longrightarrow0\,.
\end{aligned}\right.
\end{equation}
Here, we consider the following result of Damascelli--Ramaswamy~\cite{DamRam} and Li~\cite{Li}: if $1<p\le2$, $f$ is a locally Lipschitz continuous function in $\(0,\infty\)$ such that
\begin{equation}\label{Eq4}
\frac{f\(v\)-f\(u\)}{v-u}\le\Lambda\max\(u^\alpha,v^\alpha\)\,\,\,\,\forall u,v \text{ such that }0<u<v<s_0
\end{equation}
for some real numbers $\Lambda,s_0>0$, and $\alpha>p-2$, and $u$ is a $C^1$--solution of \eqref{Eq3} such that 
\begin{align}
&u\(x\)=\OO\(\left|x\right|^{-m}\)\quad\text{and}\quad\left|\nabla u\(x\)\right|=\OO\(\left|x\right|^{-m-1}\)\label{Eq5}\\
&(\text{and }u\(x\)\ge C\left|x\right|^{-m}\text{ for large }\left|x\right|\text{ when }\alpha<0)\label{Eq6}
\end{align}
for some real numbers $C>0$ and $m>p/\(\alpha+2-p\)$, then $u$ is radially symmetric and strictly radially decreasing about some point $x_0\in\R^n$, i.e. there exists $v\in C^1\(0,\infty\)$ such that $v'\(r\)<0$ for all $r>0$ and $u\(x\)=v\(\left|x-x_0\right|\)$ for all $x\in\R^n$. We also mention that other symmetry results for problems of type \eqref{Eq3} have been established without any decay assumption in the case where $f$ is nonincreasing near 0 (see Gidas--Ni--Nirenberg~\cite{GidNiNir}, Li~\cite{Li}, and Li--Ni~\cite{LiNi} in case $p=2$, Damascelli--Pacella--Ramaswamy~\cite{DamPacRam}, Damascelli--Ramaswamy~\cite{DamRam}, and Serrin--Zou~\cite{SerZou} in case $p\ne2$).

\smallskip
In case $\alpha=p^*-2$, the conditions \eqref{Eq5}--\eqref{Eq6} follow from \eqref{ThEq1}--\eqref{ThEq2} with $m=\(n-p\)/\(p-1\)$ (which is greater than $p/\(\alpha+2-p\)=\(n-p\)/p$). Consequently, by combining Theorem~\ref{Th}, the results of Damascelli--Ramaswamy~\cite{DamRam} and Li~\cite{Li}, and the regularity results that are referred to in Lemma~\ref{Lem1} below, we obtain the following corollary.

\begin{corollary}\label{Cor1}
Assume that $1<p\le2$. Let $f$ be a locally Lipschitz continuous function in $\(0,\infty\)$ such that \eqref{Eq2} and \eqref{Eq4} hold true with $f\(x,s\)=f\(s\)$. Then any nonnegative solution of \eqref{Eq1} is radially symmetric and strictly radially decreasing about some point $x_0\in\R^n$.
\end{corollary}

Let us now comment on the positive solutions of the equation with pure power nonlinearity, namely
\begin{equation}\label{Eq7}
-\Delta_pu=u^{p^*-1},\quad u>0\quad\text{in }\R^n.
\end{equation}
Guedda--V\'eron~\cite{GueVer} proved that the only positive, radially symmetric solutions of \eqref{Eq7} are of the form
\begin{equation}\label{Eq8}
u_{\mu,x_0}\(x\)=\(n\mu\)^{\frac{n-p}{p^2}}\(\frac{n-p}{p-1}\)^{\frac{\(n-p\)\(p-1\)}{p^2}}\(\mu+\left|x-x_0\right|^{\frac{p}{p-1}}\)^{\frac{p-n}{p}}
\end{equation}
for all $x\in\R^n$, for some real number $\mu>0$ and point $x_0\in\R^n$. In case $p=2$, Caffarelli--Gidas--Spruck~\cite{CafGidSpr} (see also Chen--Li~\cite{ChenLi}) proved that the functions \eqref{Eq8} are the only positive solutions of \eqref{Eq7}. In a recent paper, Damascelli--Merch\'an--Montoro--Sciunzi~\cite{DamMerMonSci} proved that any solution in $D^{1,p}\(\R^n\)$ of \eqref{Eq7} is radially symmetric provided that $2n/\(n+2\)\le p<2$. The condition $p\ge2n/\(n+2\)$ corresponds to the values of $p$ for which the function $s\mapsto s^{p^*-1}$ is Lipschitz continuous near 0. With the above Corollary~\ref{Cor1}, we extend the result of Damascelli--Merch\'an--Montoro--Sciunzi~\cite{DamMerMonSci} to the whole interval $1<p<2$. By combining the result of Guedda--V\'eron~\cite{GueVer} and Corollary~\ref{Cor1}, we obtain the following corollary.

\begin{corollary}\label{Cor2}
Assume that $1<p<2$. Then the functions \eqref{Eq8} are the only positive solutions in $D^{1,p}\(\R^n\)$ of \eqref{Eq7}. 
\end{corollary}

As a motivation to our results, it is well known that the profile of solutions of the equation
\begin{equation}\label{Eq9}
-\Delta_pu=\left|u\right|^{p^*-2}u\quad\text{in }\R^n
\end{equation}
plays a central role in the blow-up theories of critical equations. Possible references in book form on this subject and its applications in case $p=2$ are Druet--Hebey--Robert~\cite{DruHebRob}, Ghoussoub~\cite{Gho}, and Struwe~\cite{Str2}. In case $p\ne2$, global compactness results in energy spaces in the spirit of Struwe~\cite{Str1} have been established in different contexts by Alves~\cite{Alv} for equations posed in the whole $\R^n$, Saintier~\cite{Sai} in the case of a smooth, compact manifold, and Mercuri--Willem~\cite{MerWil} and Yan~\cite{Yan} in the case of a smooth, bounded domain. In view of these results, it is likely that the new information provided by Theorem~\ref{Th} and Corollary~\ref{Cor1} on the solutions of \eqref{Eq9} will lead to new existence and multiplicity results as it is the case for $p=2$.

\smallskip
The paper is organized as follows. Section~\ref{Sec2} is concerned with global boundedness results. The key result in this section is a global bound in weak Lebesgue spaces which we obtain by arguments of measure theory. In Section~\ref{Sec3}, we establish a preliminary decay estimate which is not sharp but which turns out to be a crucial ingredient in what follows. To prove this estimate, we exploit the rescaling law of the equation, and we apply a doubling property of Pol\'a\v{c}ik--Quittner--Souplet~\cite{PolQuiSou}. In Section~\ref{Sec4}, we conclude the proof of Theorem~\ref{Th}. The proof of the upper bound estimates \eqref{ThEq1} follows from the results of Sections~\ref{Sec2} and~\ref{Sec3} together with Harnack-type inequalities of Serrin~\cite{Ser} and Trudinger~\cite{Tru2}. The proof of the lower bound estimate \eqref{ThEq2} requires a Harnack inequality on annuli in the spirit of Friedman--V{\'e}ron~\cite{FriVer}.

\medskip\noindent
{\bf Acknowledgments.} The author wishes to express his gratitude to Emmanuel Hebey for helpful comments on the manuscript.

\section{Global boundedness results}\label{Sec2}

The first result of this section refers to some known regularity results for critical equations.

\begin{lemma}\label{Lem1}
Let $f:\R^n\times\R\to\R$ be a Caratheodory function such that \eqref{Eq2} holds true. Then any solution of \eqref{Eq1} belongs to $W^{1,\infty}\(\R^n\)\cap C^{1,\theta}_{\loc}\(\R^n\)$ for some $\theta\in\(0,1\)$.
\end{lemma}

\proof[Proof of Lemma~\ref{Lem1}]
A straightforward adaptation of Peral~\cite{Per}*{Theorem~E.0.20} (which in turn is adapted from Trudinger~\cite{Tru2}*{Theorem~3}) yields that for any solution $u$ of \eqref{Eq1}, there exist constants $C,R>0$ and  $\beta>1$ such that $\left\|u\right\|_{L^{\beta p^*}\(B\(x,R\)\)}\le C$ for all $x\in\R^n$, where $B\(x,R\)$ is the Euclidean ball of center $x$ and radius $R$. We then obtain a global $L^\infty$--bound by applying Serrin~\cite{Ser}*{Theorem~1}.

Once we have the $L^\infty$--boundedness of the solutions, the results of DiBenedetto~\cite{DiB} and Tolksdorf~\cite{Tol} provide global $L^\infty$--bounds and local H\"older regularity for the derivatives.
\endproof

The next result is concerned with the boundedness of solutions of \eqref{Eq1} in weak Lebesgue spaces. For any $s\in\(0,\infty\)$ and any domain $\Omega\subset\R^n$, we define $L^{s,\infty}\(\Omega\)$ as the set of all measurable functions $u:\Omega\to\R$ such that
$$\left\|u\right\|_{L^{s,\infty}\(\Omega\)}:=\sup_{h>0}\big(h\cdot\meas\(\left\{\left|u\right|>h\right\}\)^{1/s}\big)<\infty\,,$$
where $\meas\(\left\{\left|u\right|>h\right\}\)$ is the measure of the set $\left\{x\in\Omega:\,\left|u\(x\)\right|>h\right\}$. The map $\left\|\cdot\right\|_{L^{s,\infty}\(\Omega\)}$ defines a quasi-norm on $L^{s,\infty}\(\Omega\)$ (see for instance Grafakos~\cite{Gra}). 

\smallskip
Our result is as follows.

\begin{lemma}\label{Lem2}
Let $f:\R^n\times\R\to\R$ be a Caratheodory function such that \eqref{Eq2} holds true. Then any solution of \eqref{Eq1} belongs to $L^{p_*-1,\infty}\(\R^n\)$, where $p_*:=p\(n-1\)/\(n-p\)$, and hence by interpolation, any solution of \eqref{Eq1} belongs to $L^s\(\R^n\)$ for all $s\in\(p_*-1,\infty\]$.
\end{lemma}

\proof[Proof of Lemma~\ref{Lem2}]
We let $u$ be a nontrivial solution of \eqref{Eq1}. For any $h>0$, by testing \eqref{Eq1} with $T_h\(u\):=\sgn\(u\)\cdot\min\(\left|u\right|,h\)$, where $\sgn\(u\)$ denotes the sign of $u$, we obtain
\begin{equation}\label{Lem2Eq1}
\int_{\left|u\right|\le h}\left|\nabla u\right|^pdx=\int_{\left|u\right|\le h}f\(x,u\)\cdot u\,dx+h\int_{\left|u\right|>h}f\(x,u\)\cdot\sgn\(u\)dx\,.
\end{equation}
It follows from \eqref{Eq2} and \eqref{Lem2Eq1}that
\begin{equation}\label{Lem2Eq2}
\int_{\left|u\right|\le h}\left|\nabla u\right|^pdx\le\Lambda\(\int_{\left|u\right|\le h}\left|u\right|^{p^*}dx+h\int_{\left|u\right|>h}\left|u\right|^{p^*-1}dx\).
\end{equation}
We then write 
\begin{equation}\label{Lem2Eq3}
\int_{\left|u\right|\le h}\left|u\right|^{p^*}dx=\int_{\R^n}\left|T_h\(u\)\right|^{p^*}dx-h^{p^*}\meas\(\left\{\left|u\right|>h\right\}\)
\end{equation}
and
\begin{multline}\label{Lem2Eq4}
\int_{\left|u\right|>h}\left|u\right|^{p^*-1}dx=\(p^*-1\)\int_0^\infty s^{p^*-2}\meas\(\left\{\left|u\right|>\max\(s,h\)\right\}\)ds\\
=h^{p^*-1}\meas\(\left\{\left|u\right|>h\right\}\)+\(p^*-1\)\int_h^\infty s^{p^*-2}\meas\(\left\{\left|u\right|>s\right\}\)ds\,.
\end{multline}
It follows from \eqref{Lem2Eq2}--\eqref{Lem2Eq4} that
\begin{multline}\label{Lem2Eq5}
\int_{\left|u\right|\le h}\left|\nabla u\right|^pdx\le\Lambda\bigg(\int_{\R^n}\left|T_h\(u\)\right|^{p^*}dx\\
+\(p^*-1\)h\int_h^\infty s^{p^*-2}\meas\(\left\{\left|u\right|>s\right\}\)ds\bigg).
\end{multline}
Sobolev inequality gives
\begin{equation}\label{Lem2Eq6}
\int_{\R^n}\left|T_h\(u\)\right|^{p^*}dx\le K\(\int_{\left|u\right|\le h}\left|\nabla u\right|^pdx\)^{\frac{n}{n-p}}
\end{equation}
for some constant $K=K\(n,p\)$. By \eqref{Lem2Eq3}, \eqref{Lem2Eq5}, \eqref{Lem2Eq6}, and since $\int_{\R^n}\left|T_h\(u\)\right|^{p^*}dx=\oo\(1\)$ as $h\to0$, we obtain
\begin{align}\label{Lem2Eq7}
h^{p^*}\meas\(\left\{\left|u\right|>h\right\}\)&\le\int_{\R^n}\left|T_h\(u\)\right|^{p^*}dx\nonumber\\
&\le C\(h\int_h^\infty s^{p^*-2}\meas\(\left\{\left|u\right|>s\right\}\)ds\)^{\frac{n}{n-p}}
\end{align}
for small $h$, for some constant $C=C\(n,p,\Lambda\)$. We then define
$$F\(h\):=\(\int_h^\infty f\(s\)ds\)^{\frac{-p}{n-p}},\quad\text{where }f\(s\):= s^{p^*-2}\meas\(\left\{\left|u\right|>s\right\}\).$$
Since the function $t\mapsto t^{-p/\(n-p\)}$ is locally Lipschitz in $\(0,\infty\)$ and $\int_h^\infty f\(s\)ds>0$ for all $h<\left\|u\right\|_{L^\infty\(\R^n\)}$, we get that $F$ is locally absolutely continuous in $\big(0,\left\|u\right\|_{L^\infty\(\R^n\)}\big)$ with derivative
\begin{equation}\label{Lem2Eq8}
F'\(h\)=\frac{p}{n-p}\(\int_h^\infty f\(s\)ds\)^{\frac{-n}{n-p}}f\(h\)
\end{equation}
for a.e. $h\in\big(0,\left\|u\right\|_{L^\infty\(\R^n\)}\big)$ (see for instance Leoni~\cite{Leoni}*{Theorem 3.68}). By \eqref{Lem2Eq7} and \eqref{Lem2Eq8}, we obtain
\begin{equation}\label{Lem2Eq9}
F'\(h\)\le C\cdot\frac{p}{n-p}\cdot h^{\frac{2p-n}{n-p}}
\end{equation}
for small $h$. Integrating \eqref{Lem2Eq9} gives
\begin{equation}\label{Lem2Eq10}
F\(h\)-F\(0\)\le Ch^{\frac{p}{n-p}}
\end{equation}
for small $h$, where $F\(0\):=\lim_{h\to0}F\(h\)$. On the other hand, by \eqref{Lem2Eq4} and dominated convergence, we have
\begin{equation}\label{Lem2Eq11}
\(p^*-1\)hF\(h\)^{\frac{p-n}{p}}\le h\int_{\left|u\right|>h}\left|u\right|^{p^*-1}dx=\oo\(1\)
\end{equation}
as $h\to0$. It follows from \eqref{Lem2Eq10} and \eqref{Lem2Eq11} that $F\(0\)>0$, i.e. $\int_0^\infty f\(s\)ds<\infty$. By \eqref{Lem2Eq7} and since $p^*-\frac{n}{n-p}=p_*-1$ and $F$ is nonincreasing, we then get 
$$h^{p_*-1}\meas\(\left\{\left|u\right|>h\right\}\)\le CF\(h\)^{-n/p}\le CF\(0\)^{-n/p}$$
for small $h$, and hence we obtain $\left\|u\right\|_{L^{p_*-1,\infty}\(\R^n\)}<\infty$.
\endproof

By \eqref{Eq2} and a weak version of Kato's inequality~\cite{Kato} (see Leon~\cite{Leon}*{Proposition~3.2}), we obtain
\begin{equation}\label{Eq10}
-\Delta_p\left|u\right|\le\left|f\(x,u\)\right|\le\Lambda\left|u\right|^{p^*-1}\quad\text{in }\R^n,
\end{equation}
where the inequality is in the sense that
$$\int_{\R^n}\left|\nabla\left|u\right|\right|^{p-2}\nabla\left|u\right|\cdot\nabla\varphi\,dx\le\Lambda\int_{\R^n}\left|u\right|^{p^*-1}\varphi\,dx$$
for all nonnegative, smooth function $\varphi$ with compact support in $\R^n$.

\smallskip
Our last result in this section is as follows.

\begin{lemma}\label{Lem3}
For any real number $\Lambda>0$ and any nonnegative, nontrivial solution $v\in D^{1,p}\(\R^n\)$ of the inequality $-\Delta_pv\le\Lambda v^{p^*-1}$ in $\R^n$, we have $\left\|v\right\|_{L^{p^*}\(\R^n\)}\ge\kappa_0$ for some constant $\kappa_0=\kappa_0\(n,p,\Lambda\)>0$.
\end{lemma}

\proof
By Sobolev inequality and testing the inequality $-\Delta_pv\le\Lambda v^{p^*-1}$ with the function $v$, we obtain
\begin{equation}\label{Lem3Eq}
\int_{\R^n}v^{p^*}dx\le K\(\int_{\R^n}\left|\nabla v\right|^pdx\)^{\frac{n}{n-p}}\le K\(\Lambda\int_{\R^n}v^{p^*}dx\)^{\frac{n}{n-p}}
\end{equation}
for some constant $K=K\(n,p\)$. The result then follows immediately from \eqref{Lem3Eq}.
\endproof

\section{A preliminary decay estimate}\label{Sec3}

The following result provides a decay estimate which is not sharp but which will serve as a preliminary step in the proof of Theorem~\ref{Th}. The proof relies on rescaling arguments and a doubling property of Pol\'a\v{c}ik--Quittner--Souplet~\cite{PolQuiSou}.

\begin{lemma}\label{Lem4}
Let $\kappa_0$ be as in Lemma~\ref{Lem3}, $f:\R^n\times\R\to\R$ be a Caratheodory function such that \eqref{Eq2} holds true, and $u$ be a solution of \eqref{Eq1}. For any $\kappa>0$, we define
$$r_\kappa\(u\):=\inf\big(\big\{r>0\,:\,\left\|u\right\|_{L^{p^*}\(\R^n\backslash B\(0,r\)\)}<\kappa\big\}\big),$$
where $B\(0,r\)$ denotes the Euclidean ball of center 0 and radius $r$. Then for any $\kappa\in\(0,\kappa_0\)$ and $r>r_\kappa\(u\)$, there exists a constant $K_0=K_0\(n,p,\Lambda,\kappa,r,r_\kappa\(u\)\)$ such that
\begin{equation}\label{Lem4Eq1}
\left|u\(x\)\right|\le K_0\left|x\right|^{\frac{p-n}{p}}\quad\text{for all }x\in\R^n\backslash B\(0,r\).
\end{equation}
\end{lemma}

\proof[Proof of Lemma~\ref{Lem4}]
We fix $\Lambda>0$, $\kappa\in\(0,\kappa_0\)$, $r>0$, and $r'\in\(0,r\)$. As is easily seen, in order to prove Lemma~\ref{Lem4}, it is sufficient to prove that there exists a constant $K_1=K_1\(n,p,\Lambda,\kappa,r,r'\)$ such that for any solution $u$ of \eqref{Eq1} such that $r_\kappa\(u\)\le r'$, we have
\begin{equation}\label{Lem4Eq2}
\dist\(x,B\(0,r''\)\)\left|u\(x\)\right|^{\frac{p}{n-p}}\le K_1\quad\text{for all }x\in\R^n\backslash B\(0,r\),
\end{equation}
where $r'':=\(r+r'\)/2$ and $\dist$ is the Euclidean distance function. 

We prove \eqref{Lem4Eq2} by contradiction. Suppose that for any $\alpha\in\N$, there exists a Caratheodory function $f_\alpha:\R^n\times\R\to\R$ such that \eqref{Eq2} holds true, a solution $u_\alpha$ of \eqref{Eq1} with $f=f_\alpha$ such that $r_\kappa\(u_\alpha\)\le r'$, and a point $x_\alpha\in\R^n\backslash B\(0,r\)$ such that
\begin{equation}\label{Lem4Eq3}
\dist\(x_\alpha,B\(0,r''\)\)\left|u_\alpha\(x_\alpha\)\right|^{\frac{p}{n-p}}>2\alpha\,.
\end{equation}
By \eqref{Lem4Eq3} and Pol\'a\v{c}ik--Quittner--Souplet~\cite{PolQuiSou}*{Lemma~5.1}, we get that there exists $y_\alpha\in\R^n\backslash B\(0,r''\)$ such that
\begin{equation}\label{Lem4Eq4}
\dist\(y_\alpha,B\(0,r''\)\)\left|u_\alpha\(y_\alpha\)\right|^{\frac{p}{n-p}}>2\alpha\,,\quad\left|u_\alpha\(x_\alpha\)\right|\le\left|u_\alpha\(y_\alpha\)\right|,
\end{equation}
and 
\begin{equation}\label{Lem4Eq5}
\left|u_\alpha\(y\)\right|\le2^{\frac{n-p}{p}}\left|u_\alpha\(y_\alpha\)\right|\quad\text{for all }y\in B\big(y_\alpha,\alpha\,\left|u_\alpha\(y_\alpha\)\right|^{\frac{-p}{n-p}}\big).
\end{equation}
For any $\alpha$ and $y\in\R^n$, we define
\begin{equation}\label{Lem4Eq6}
\widetilde{u}_\alpha\(y\):=\mu_\alpha\cdot u_\alpha\big(\mu_\alpha^{\frac{p}{n-p}}\cdot y+y_\alpha\big),
\end{equation}
where $\mu_\alpha:=\left|u_\alpha\(y_\alpha\)\right|^{-1}$. By rescaling \eqref{Eq1}, we obtain
\begin{equation}\label{Lem4Eq7}
-\Delta_p\widetilde{u}_\alpha=\mu_\alpha^{p^*-1}\cdot f_\alpha\big(\mu_\alpha^{\frac{p}{n-p}}\cdot y+y_\alpha,\mu_\alpha^{-1}\cdot\widetilde{u}_\alpha\big)\quad\text{in }\R^n.
\end{equation}
It follows from \eqref{Eq2} that
\begin{equation}\label{Lem4Eq8}
\big|\mu_\alpha^{p^*-1}\cdot f_\alpha\big(\mu_\alpha^{\frac{p}{n-p}}\cdot y+y_\alpha,\mu_\alpha^{-1}\cdot\widetilde{u}_\alpha\big)\big|\le\Lambda\left|\widetilde{u}_\alpha\right|^{p^*-1}\quad\text{in }\R^n.
\end{equation}
Moreover, by \eqref{Lem4Eq5} and \eqref{Lem4Eq6}, we obtain
\begin{equation}\label{Lem4Eq9}
\left|\widetilde{u}_\alpha\(0\)\right|=1\quad\text{and}\quad\left|\widetilde{u}_\alpha\(y\)\right|\le2^{\frac{n-p}{p}}\quad\text{for all }y\in B\(0,\alpha\).
\end{equation} 
By DiBenedetto~\cite{DiB} and Tolksdorf~\cite{Tol}, it follows from \eqref{Lem4Eq8} and \eqref{Lem4Eq9} that there exists a constant $C>0$ and a real number $\theta\in\(0,1\)$ such that for point $x\in\R^n$, we have
\begin{equation}\label{Lem4Eq10}
\left\|\widetilde{u}_\alpha\right\|_{C^{1,\theta}(B\(x,1\))}\le C
\end{equation}
for large $\alpha$. By compactness of $C^{1,\theta}(B\(x,1\))\hookrightarrow C^1(B\(x,1\))$, it follows from \eqref{Lem4Eq10} that $\(\widetilde{u}_\alpha\)_\alpha$ converges up to a subsequence in $C^1_{\loc}\(\R^n\)$ to some function $\widetilde{u}_\infty$. By \eqref{Lem4Eq9}, we obtain $\left|\widetilde{u}_\infty\(0\)\right|=1$. Moreover, by observing that both the $D^{1,p}$--norm and the inequality \eqref{Eq10} are left invariant by the rescaling \eqref{Lem4Eq7}, we get that $\widetilde{u}_\infty\in D^{1,p}\(\R^n\)$ and $\left|\widetilde{u}_\infty\right|$ is a weak solution of
\begin{equation}\label{Lem4Eq11}
-\Delta_p\left|\widetilde{u}_\infty\right|\le\Lambda\left|\widetilde{u}_\infty\right|^{p^*-1}\quad\text{in }\R^n.
\end{equation}

On the other hand, for any $R>0$, we have
\begin{equation}\label{Lem4Eq12}
\left\|\widetilde{u}_\alpha\right\|_{L^{p^*}(B\(0,R\))}=\left\|u_\alpha\right\|_{L^{p^*}(B(y_\alpha,R\mu_\alpha^{\frac{p}{n-p}}))}.
\end{equation}
By \eqref{Lem4Eq4} and since $r_\kappa\(u_\alpha\)<r''$, we get
\begin{equation}\label{Lem4Eq13}
B\big(y_\alpha,R\mu_\alpha^{\frac{p}{n-p}}\big)\cap B\(0,r_\kappa\(u_\alpha\)\)=\emptyset
\end{equation}
for large $\alpha$. By \eqref{Lem4Eq12}, \eqref{Lem4Eq13}, and by definition of $r_\kappa\(u_\alpha\)$, we obtain
\begin{equation}\label{Lem4Eq14}
\left\|\widetilde{u}_\alpha\right\|_{L^{p^*}(B\(0,R\))}\le\kappa
\end{equation}
for large $\alpha$. Passing to the limit into \eqref{Lem4Eq14} as $\alpha\to\infty$ and then as $R\to\infty$ yields
\begin{equation}\label{Lem4Eq15}
\left\|\widetilde{u}_\infty\right\|_{L^{p^*}\(\R^n\)}\le\kappa\,.
\end{equation}
Since $\kappa<\kappa_0$, by Lemma~\ref{Lem3}, \eqref{Lem4Eq11}, and \eqref{Lem4Eq15}, we get that $\widetilde{u}_\infty\equiv0$, which is in contradiction with $\left|\widetilde{u}_\infty\(0\)\right|=1$. This ends the proof of Lemma~\ref{Lem4}.
\endproof

\section{Proof of Theorem~\ref{Th}}\label{Sec4}

We can now prove Theorem~\ref{Th} by applying Lemmas~\ref{Lem2},~\ref{Lem4}, and Harnack-type inequalities of Serrin~\cite{Ser} and Trudinger~\cite{Tru2}.

\proof[Proof of \eqref{ThEq1}]
We let $u$ be a solution of \eqref{Eq1}. We let $\kappa$ and $r$ be as in Lemma~\ref{Lem4}. For any $R>0$ and $y\in\R^n$, we define 
\begin{equation}\label{ThEq3}
u_R\(y\):=R^{\frac{n-p}{p-1}}\cdot u\big(R\cdot y\big).
\end{equation}
By rescaling \eqref{Eq1}, we obtain
\begin{equation}\label{ThEq4}
-\Delta_pu_R=R^n\cdot f\big(R\cdot y,R^{\frac{p-n}{p-1}}\cdot u_R\big)\quad\text{in }\R^n.
\end{equation}
It follows from \eqref{Eq2} that
\begin{equation}\label{ThEq5}
\big|R^n\cdot f\big(R\cdot y,R^{\frac{p-n}{p-1}}\cdot u_R\big)\big|\le\Lambda\cdot R^{\frac{-p}{p-1}}\cdot\left|u_R\right|^{p^*-1}\quad\text{in }\R^n.
\end{equation}
Moreover, similarly to \eqref{Eq10}, it follows from \eqref{ThEq4} and \eqref{ThEq5} that $\left|u_R\right|$ is a weak solution of
\begin{equation}\label{ThEq6}
-\Delta_p\left|u_R\right|\le\Lambda\cdot R^{\frac{-p}{p-1}}\cdot\left|u_R\right|^{p^*-1}\quad\text{in }\R^n.
\end{equation}
By writing $\left|u_R\right|^{p^*-1}=\left|u_R\right|^{p^*-p}\cdot\left|u_R\right|^{p-1}$ and applying Lemma~\ref{Lem4}, we obtain
\begin{equation}\label{ThEq7}
R^{\frac{-p}{p-1}}\cdot\left|u_R\right|^{p^*-1}\le K_0^{p^*-p}\left|u_R\right|^{p-1}\quad\text{in }\R^n\backslash B\(0,1\)
\end{equation}
provided that $R\ge r$. It follows from \eqref{ThEq6}, \eqref{ThEq7}, and Trudinger~\cite{Tru2}*{Theorem~1.3} that for any $\varepsilon>0$, we have
\begin{equation}\label{ThEq8}
\left\|u_R\right\|_{L^\infty\(B\(0,2\)\backslash B\(0,4\)\)}\le c_\varepsilon\left\|u_R\right\|_{L^{p-1+\varepsilon}\(B\(0,5\)\backslash B\(0,1\)\)}.
\end{equation}
for some constant $c_\varepsilon=c\(n,p,\Lambda,K_0,\varepsilon\)$. We fix $\varepsilon_0=\varepsilon_0\(n,p\)$ such that $0<\varepsilon_0<p_*-p$, where $p_*$ is as in Lemma~\ref{Lem2}. By a generalized version of H\"older's inequality (see for instance Grafakos~\cite{Gra}*{Exercise~1.1.11}), we obtain that there exists a constant $c_0=c_0\(n,p\)$ such that
\begin{equation}\label{ThEq9}
\left\|u_R\right\|_{L^{p-1+\varepsilon_0}\(B\(0,5\)\backslash B\(0,1\)\)}\le c_0\left\|u_R\right\|_{L^{p_*-1,\infty}\(B\(0,5\)\backslash B\(0,1\)\)}.
\end{equation}
By observing that the quasi-norm $\left\|\cdot\right\|_{L^{p_*-1,\infty}\(R^n\)}$ is left invariant by the rescaling \eqref{ThEq3}, we deduce from \eqref{ThEq8}, \eqref{ThEq9}, and Lemma~\ref{Lem2} that
\begin{equation}\label{ThEq10}
\left\|u_R\right\|_{L^\infty\(B\(0,2\)\backslash B\(0,4\)\)}\le c_1
\end{equation}
for some constant $c_1=c_1\big(n,p,\Lambda,K_0,\left\|u\right\|_{L^{p_*-1,\infty}\(\R^n\)}\big)$. By \eqref{ThEq4}--\eqref{ThEq7}, \eqref{ThEq10}, and the estimates of DiBenedetto~\cite{DiB} and Tolksdorf~\cite{Tol}, we get
\begin{equation}\label{ThEq11}
\left\|\nabla u_R\right\|_{L^\infty\(B\(0,5/2\)\backslash B\(0,7/2\)\)}\le c_2\,.
\end{equation}
for some constant $c_2=c_2\big(n,p,\Lambda,K_0,\left\|u\right\|_{L^{p_*-1,\infty}\(\R^n\)}\big)$. Finally, for any $x\in\R^n\backslash B\(0,3r\)$, by applying \eqref{ThEq10} and \eqref{ThEq11} with $R=\left|x\right|/3$, we obtain
\begin{equation}\label{ThEq12}
\left|u\(x\)\right|\le c_3\left|x\right|^{\frac{p-n}{p-1}}\quad\text{and}\quad\left|\nabla u\(x\)\right|\le c_3\left|x\right|^{\frac{1-n}{p-1}}
\end{equation}
for some constant $c_3=c_3\big(n,p,\Lambda,K_0,\left\|u\right\|_{L^{p_*-1,\infty}\(\R^n\)}\big)$. Since on the other hand $u$ and $\nabla u$ are uniformly bounded in $B\(0,3r\)$, we can deduce \eqref{ThEq1} from \eqref{ThEq12}.
\endproof

\begin{remark}\label{Rem1}
As one can see from the above proof, the constant $C_0$ in \eqref{ThEq1} depends only on $n$, $p$, $\Lambda$, $\kappa$, $r$, and upper bounds of $r_\kappa\(u\)$, $\left\|u\right\|_{L^{p_*-1,\infty}\(\R^n\)}$, and $\left\|u\right\|_{W^{1,\infty}\(B\(0,3r\)\)}$.
\end{remark}

In order to prove the lower bound estimate \eqref{ThEq2}, we need the following Harnack inequality on annuli. This result is inspired from Friedman--V{\'e}ron~\cite{FriVer} where a similar result is used for the study of singular solutions of $p$--Laplace equations in pointed domains. 

\begin{lemma}\label{Lem5}\samepage
Let $f:\R^n\times\R\to\R$ be a Caratheodory function such that \eqref{Eq2} holds true, $u$ be a nonnegative solution of \eqref{Eq1}, $\kappa$ and $r$ be as in Lemma~\ref{Lem4}, and $K_0$ be the constant given by Lemma~\ref{Lem4}. Then there exists a constant $c_4=c_4\(n,p,\Lambda,K_0\)$ such that 
\begin{equation}\label{Lem5Eq1}
\sup_{2R<\left|x\right|<5R}\(u\(x\)\)\le c_4\cdot\inf_{2R<\left|x\right|<5R}\(u\(x\)\)
\end{equation}
for all $R\ge r$.
\end{lemma}

\proof[Proof of Lemma~\ref{Lem5}]
For any $R>0$, we define $u_R$ as in \eqref{ThEq3}. By \eqref{ThEq4}, \eqref{ThEq5}, \eqref{ThEq7}, and Serrin~\cite{Ser}*{Theorem~5}, we obtain that there exists a constant $c=c\(n,p,\Lambda,K_0\)$ such that
\begin{equation}\label{Lem5Eq2}
\sup_{z\in B\(y,1/3\)}\(u_R\(z\)\)\le c\cdot\inf_{z\in B\(y,1/3\)}\(u_R\(z\)\)
\end{equation}
for all points $y$ in the annulus $A:=B\(0,5\)\backslash B\(0,2\)$. Moreover, we can join every two points in $A$ by 17 connected balls of radius 1/3 and centers in $A$. Hence \eqref{Lem5Eq1} follows from \eqref{Lem5Eq2} with $c_4:=c^{17}$.
\endproof

We can now prove \eqref{ThEq2} by applying Lemma~\ref{Lem5}.

\proof[Proof of \eqref{ThEq2}]
We let $u$ be a nonnegative solution of \eqref{Eq1} such that $\int_{\R^n}f\(x,u\)dx>0$. In particular, in view of \eqref{Eq2}, we have $u\not\equiv0$, and hence $u>0$ in $\R^n$ by the strong maximum principle of Vazquez~\cite{Vaz}. By Lemma~\ref{Lem5}, we then get that in order to prove \eqref{ThEq2}, it is sufficient to prove a lower bound estimate of $\left\|u\right\|_{L^{\infty}\(B\(0,5R\)\backslash B\(0,2R\)\)}$ for large $R$. 

By \eqref{ThEq4}, \eqref{ThEq5}, \eqref{ThEq7}, and Serrin~\cite{Ser}*{Theorem~1}, we obtain
\begin{align}\label{ThEq13}
\left\|\nabla u_R\right\|_{L^p\(B\(0,4\)\backslash B\(0,3\)\)}&\le c_5\left\|u_R\right\|_{L^p\(B\(0,5\)\backslash B\(0,2\)\)}\nonumber\\
&\le c'_5\left\|u_R\right\|_{L^{\infty}\(B\(0,5\)\backslash B\(0,2\)\)}
\end{align}
for some constants $c_5$ and $c'_5$ depending only on $n$, $p$, $\Lambda$, and $K_0$, where $u_R$ is as in \eqref{ThEq3}. Rescaling \eqref{ThEq13} yields
\begin{equation}\label{ThEq14}
\left\|\nabla u\right\|_{L^p\(B\(0,4R\)\backslash B\(0,3R\)\)}\le c'_5R^{\frac{n-p}{p}}\left\|u\right\|_{L^{\infty}\(B\(0,5R\)\backslash B\(0,2R\)\)}.
\end{equation}

Next, we claim that if $\int_{\R^n}f\(x,u\)dx>\lambda$ for some real number $\lambda>0$, then we have
\begin{equation}\label{ThEq15}
\left\|\nabla u\right\|_{L^p\(B\(0,4R\)\backslash B\(0,3R\)\)}\ge c_6R^\frac{p-n}{p\(p-1\)}
\end{equation}
for large $R$, for some constant $c_6=c_6\(n,p,\lambda\)>0$. For any $x\in\R^n$ and $R>0$, we define $\chi_R\(x\):=\chi\(\left|x\right|/R\)$, where $\chi\in C^1\(0,\infty\)$ is a cutoff function such that $\chi\equiv1$ on $\[0,3\]$, $\chi\equiv0$ on $\[4,\infty\)$, $0\le\eta\le1$ and $\left|\eta'\right|\le2$ on $\(3,4\)$. By testing \eqref{Eq1} with $\chi_R$ and applying H\"older's inequality, we obtain
\begin{align}\label{ThEq16}
\int_{\R^n}f\(x,u\)\chi_R\,dx&=\int_{\R^n}\left|\nabla u\right|^{p-2}\nabla u\cdot\nabla\chi_R\,dx\nonumber\\
&\le\left\|\nabla u\right\|^{p-1}_{L^p\(\supp\(\nabla\chi_R\)\)}
\cdot\left\|\nabla\chi_R\right\|_{L^p\(\supp\(\nabla\chi_R\)\)},
\end{align}
where $\supp\(\chi_R\)$ denotes the support of $\chi_R$. It follows from \eqref{ThEq16} and the definition of $\chi_R$ that
\begin{equation}\label{ThEq17}
\int_{\R^n}f\(x,u\)\chi_R\,dx\le CR^{\frac{n-p}{p}}\left\|\nabla u\right\|^{p-1}_{L^p\(B\(0,4R\)\backslash B\(0,3R\)\)}
\end{equation}
for some constant $C=C\(n,p\)>0$. Then \eqref{ThEq15} follows from \eqref{ThEq17} with $c_6:=\(\lambda/C\)^{\frac{1}{p-1}}$

Finally, we deduce \eqref{ThEq2} from \eqref{Lem5Eq1}, \eqref{ThEq14}, and \eqref{ThEq15}.
\endproof

\begin{remark}\label{Rem2}
As one can see from the above proof, the constant $C_1$ in \eqref{ThEq2} depends only on $n$, $p$, $\lambda$, $\Lambda$, $\kappa$, $r$, $r_\kappa\(u\)$, and a lower bound of $u$ on the ball $B\(0,2\max\(r,R_{\lambda,f}\(u\)\)\)$, where 
$$R_{\lambda,f}\(u\):=\inf\big(\big\{R>0:\,\int_{\R^n}f\(x,u\)\chi_{R'}dx>\lambda\,,\quad\forall R'>R\big\}\big).$$
\end{remark}

\end{document}